\pdfoutput=1
\documentclass{article}
\usepackage[a4paper,  top = 1.2in, bottom = 1.2in, left = 1.5in, right = 1.5in]{geometry}

\usepackage{caption}

\usepackage{amsmath,amssymb}

\usepackage{latexsym,fancybox,theorem}
\usepackage{here}
\usepackage{longtable}
\usepackage{lscape}
\usepackage{setspace}
\usepackage[dviout]{graphicx}

\usepackage{hyperref}
\usepackage{optidef}

\newtheorem{prf}{Proof.}
\def\qed{\hfill $\Box$}

  \makeatletter
  \@addtoreset{equation}{section}
  \makeatother

  \makeatletter
    
    \@addtoreset{figure}{section}
  \makeatother

 \makeatletter
    
   \@addtoreset{table}{section}
  \makeatother

\theoremstyle{plain}
\newtheorem{thm}{Theorem}[section]
\newtheorem{lem}{Lemma}[section]

\newtheorem{alg}{Algorithm}[section]
\newtheorem{asm}{Assumption}[section]

\newtheorem{definition}{Definition}[section]
\numberwithin{equation}{section}
\newcommand{\R}{\mathbb{R}}

\newcommand{\g}{\nabla}
\newcommand{\gk}{\nabla f(x_k)}

\newcommand{\h}{\nabla^2}

\newcommand{\ak}{\alpha_k}
\newcommand{\lb}{L-BFGS\ }
\newcommand{\dk}{d_k(\mu)}

\newcommand{\yk}{\hat{y}_{k}(\mu)}

\newcommand{\hk}{\hat{H}_{k}(\mu)}

\title{A Regularized Limited Memory BFGS method for Large-Scale Unconstrained Optimization 
and its efficient Implementations}

\author{Hardik TANKARIA \and Shinji SUGIMOTO \and 
Nobuo YAMASHITA\footnote{Department of Applied Mathematics and Physics, Graduate 
School of Informatics, Kyoto University, Yoshida-Honmachi, Sakyo-ku, Kyoto 606-8501, Japan. E-mail: \tt{nobuo@i.kyoto-u.ac.jp}}}

\begin{document}
\newpage
\maketitle


\begin{abstract}

The limited memory BFGS~(L-BFGS)  method is one of the popular methods for solving
large-scale unconstrained optimization. Since the standard L-BFGS method uses a line search to 
guarantee its global convergence, it sometimes requires a large number of function 
evaluations. To overcome the difficulty, we propose a new L-BFGS with a certain 
regularization technique. We show its global convergence  under the usual assumptions. In order to make 
the method more robust and efficient, we also extend it with several techniques such as 
nonmonotone technique and simultaneous use of the Wolfe line search. Finally, we present 
some numerical results for test problems in CUTEst,  which show that the proposed method is robust in terms 
of solving more number of problems.
\end{abstract}

\section{Introduction}
In this paper we consider the large-scale unconstrained optimization problem:


\begin{mini}|l|
{x\in\mathbb{R}^n}{f(x)}{}{},
\label{1}
\end{mini}
where $f: \R^n \to \R$ is a smooth function. For solving it, we focus on the quasi-Newton  
type method as
\begin{equation*}
 x_{k+1}=x_k + d_k,
\end{equation*}
where $x_k \in \R^n $ is the $k^{th}$ iteration and $d_k \in \R^n$ denotes a search direction obtained by 
a certain quasi-Newton method.

The standard solution methods to solve \eqref{1}  such as the steepest descent method, 
Newton's method and the BFGS method \cite{Dennis,Nocedalbook} are not suitable for large-scale problems.
 This is  because the steepest descent method generally converges slowly, while Newton's 
 method needs to  compute the Hessian matrix and solve  linear equations at each iteration. 
 Moreover, the BFGS  method requires $O(n^2) $ memory to store and calculate the approximate
 Hessian of $f$, which causes some difficulty for large-scale problem.

One of the popular quasi-Newton methods for solving large-scale problem is the limited 
memory BFGS(L-BFGS)  \cite{Liu,Nocedal}, which uses small memory to store an approximate Hessian of 
$f$. The L-BFGS method stores the last $m$ vector pairs of $(s_{k-i},y_{k-i}), $ $ i=0,1, \ldots,m-1,$ to 
compute a search direction $d_k$, where

\begin{align}
 s_k =& x_k-x_{k-1}, \text{\ and \ } y_k = \nabla f(x_k)-\nabla f(x_{k-1}) ,\nonumber
 \end{align}
and computes $d_k$ in $O(mn)$ time.

The usual L-BFGS adopts the Wolfe line search to guarantee its global convergence. The line search
 sometimes needs a large number of function evaluations. Thus, it is preferable to reduce the 
 number of function evaluations as much as possible.

The trust region method (TR-method) can guarantee the global convergence.
It is known that the TR-method needs fewer function evaluations than the line search 
\cite{Burke1,Burke2,More}. The L-BFGS method combined with the TR-method \cite{Burke1,Burke2} produces good 
performance for many benchmark problems in terms of the number of function evaluations. 
However, the TR-method must solve the constrained subproblem
\begin{align}
\label{2}
& \text{minimize} \quad f(x_k) + \nabla f(x_k)^T d + \frac{1}{2} d^T B_k d \\
& \text{subject to} \quad \|d\| \leq \Delta_k, \nonumber
\end{align}
in each step, where $\Delta_k$ is the trust-region radius and $B_k$ is an approximate Hessian obtained by 
L-BFGS. It takes a considerable amount of time to solve \eqref{2}.

To overcome the 
difficulty we consider adopting a regularization technique instead of the TR-method. This is motivated by
the  regularized Newton method proposed by Ueda and Yamashita \cite{Ueda1,Ueda2,Ueda3}. The 
method computes a search direction $d_k$ as a solution of the following linear equations:
\begin{equation}
\label{3}
 (\h f(x_k)+ \mu_k I) d= -\g f(x_k),
\end{equation}
where $\mu_k>0$  is called a regularized parameter. If $\mu_k$ coincides with the value of the optimal
Lagrange multiplier at a solution of problem \eqref{2}, then $d$ is a solution of  \eqref{2}. Note that 
the linear  equations  \eqref{3} are simpler than subproblem \eqref{2} of the TR-method. The regularized Newton 
method \cite{Ueda1} controls the parameter $\mu_k$ instead of computing the step length to guarantee global 
convergence. However, since the regularized Newton method in \cite{Ueda1} is based on  Newton's method, it 
must compute the Hessian matrix of $f$.

In this paper we propose a novel approach that combines the L-BFGS method with the 
regularization technique. We call the proposed method regularized L-BFGS method. 
One of natural ways to implement the idea is to use a solution of the following equations as a search direction,  
\begin{equation}
\label{4}
 (B_k + \mu I) d = -\gk,
\end{equation}
where $B_k$ is an approximate Hessian given by a certain quasi-Newton method. However when $B_k
$ is calculated by the L-BFGS method, it is difficult to compute $(B_k+\mu I)^{-1}$. Therefore, 
we try to  directly  
construct $(B_k+\mu I)^{-1}$ by the L-BFGS method for $f(x)+\mu \|x\|^2 $, that is we use $(s_k,\hat{y}_k(\mu))$, 
where $\hat{y}_k(\mu) = y_k + \mu s_k$, instead of $(s_k,y_k)$. Note that the term 
$\mu s_k$ in $\hat{y}_k(\mu)$ plays the role of regularization. Then, the search direction $d_k$ 
can be computed in  $O(mn)$ time like the conventional L-BFGS  method. For global 
convergence, we also control the regularized parameter $\mu_k$ in a way similar to the regularized 
Newton 
method \cite{Ueda1}. We then show that the proposed algorithm ensures global convergence. 

A drawback of the proposed method is that a step $d_k$ sometimes becomes small, and it causes a large number of 
iterations. To get a longer step, we propose two techniques: a nonmonotone technique  and a simultaneous use 
of the Wolfe line search. Recall that the step length given by the Wolfe condition is allowed to be larger than 1, 
and hence the step can explore a larger area. Thus, if $f(x_k+\alpha d_k)<f(x_k+d_k)$ for 
$\alpha>1$, it would be reasonable to find $\alpha$ via the Wolfe line search.

The paper is organized as follows. The regularized \lb\  is presented in section 2, and  
its global convergence is shown in section 3. In section 4 we discuss some 
implementation issues, such as a simultaneous use of R\lb and a 
nonmonotone technique. In section 5, we present numerical results by comparing three 
algorithms: the L-BFGS, the regularized L-BFGS, and the regularized \lb\ with line search. 
Section 6  concludes the paper.

Throughout the paper, we use the following notations.
For a vector $x \in \mathbb{R}^n$, $\|x\|$ denotes the Euclidean norm defined by $\|x\| := \sqrt{x^Tx} $. For a 
symmetric matrix $M \in \mathbb{R}^{n\times n}$, we denote the maximum and  minimum eigenvalues of $M$ as $
\lambda_{\max}(M)$ and $\lambda_{\min}(M)$.
Moreover, $\|M\|$ denotes the $l_2$ norm of $M$ defined by $\|M\| := \sqrt{\lambda_{\max}(M^TM)}$. If $M$ is a 
symmetric positive-semidefinite matrix, then $\|M\| = \lambda_{\max}(M)$.
Next, we give a definition of Lipschitz continuity.
\begin{definition}[Lipschitz continuity]
 Let $S$ be a subset of $\R^n$ and $f: S \to \R$.
\end{definition}
\begin{itemize}
 \item[i)] The function $f$ is said to be Lipschitz continuous on $S$ if there exists a 
 positive constant $L_f$ such that
 $$ |f(x) - f(y)| \leq L_f \| x- y \|\  \forall x,y \in S.$$
\item[ ii)] Suppose that the function $f$ is differentiable. $\g f$ is said to be Lipschitz
continuous on $S$ if there exists a positive constant $L_g$ such that 
$$ \|\g f(x) - \g f(y)\| \leq L_g \| x-y\|\  \forall x,y \in S.$$
\end{itemize}



\section{The regularized L-BFGS method}

In this section, we propose a regularized L-BFGS method that controls
the regularized parameter at each iteration.
In the following, $x_k$ denotes the $k$-th iterative point,
$B_k$ denotes the approximate Hessian of $f(x_k)$,
and $H_k^{-1} = B_k$.

We consider combining the L-BFGS method with the regularized Newton method (\ref{3}).
For this purpose, we may replace the Hessian $\nabla^2 f(x_k)$ in equation (\ref{3}) with the 
approximate Hessian $B_k$, that is, we define a search direction
$d_k$ as a solution of
\begin{eqnarray}
  (B_k+ \mu I)^{-1}d_k = -\nabla f(x_k) .\label{lineeq2}
\end{eqnarray}
However, since the L-BFGS method updates $H_k$,
it is not easy to construct $B_k$ explicitly.
Furthermore, even if we obtain $B_k$, it takes a considerable amount of time to solve the 
linear equation (\ref{lineeq2}) in large-scale cases.

Now, we may regard $B_k+\mu I$ as an approximation of $\nabla^2 f(x) + \mu I$.
Since $B_k$ is the approximate Hessian of $f(x_k)$,
the matrix $B_k+\mu I$ is an approximate Hessian of $f(x)+\frac{\mu}{2}\|x\|^2$.
The L-BFGS method uses the vector pair $(s_k, y_k)$ to construct the approximate Hessian, 
where $s_k = x_{k+1} - x_{k}$ and $y_k = \nabla f(x_{k+1}) - \nabla f(x_k)$.
Note that $y_k$ consists of the gradients of $f$. Therefore, when we compute the approximate 
Hessian of $f(x) + \frac{\mu}{2}\|x\|^2$, we use the gradients of $f(x)+\frac{\mu}{2}\|x\|^2$. 
That is, we adopt the following $\hat{y}_k(\mu)$ instead of $y_k$:
\begin{eqnarray}
 \hat{y}_k(\mu) &  =  (\nabla f(x_{k+1})+\mu x_{k+1} )-(\nabla f(x_k)+\mu x_{k})  =  y_k +  
 \mu  s_k. \nonumber
\end{eqnarray}
Let $\hat{H}_k(\mu)$ be a matrix constructed by  the L-BFGS method with vector pairs $(s_i, 
\hat{y}_i(\mu)),$ $ i=1,\ldots,m $ and an appropriate initial matrix $\hat{H}_k^{(0)}(\mu)$.
Then, the search direction $d_k = -\hat{H}_k(\mu)\nabla f(x_k)$ is calculated in $O(mn)$ time, 
which is the same as the original L-BFGS.

Note that if $s_k^T\hat{y}_k(\mu) > 0$
and $\hat{H}_k^{(0)}$ is positive-definite, then  $\hat{H}_k(\mu)$ is positive definite.
When $s_k^T \hat{y}_k(\mu) > 0$ is not satisfied, we may replace $\hat{y}_k(\mu)$ by 
$\tilde{y}_k(\mu)$:
\[
   \tilde{y}_k(\mu) = y_k + \left( \max \left\{ 0,\frac{-s_k^Ty_k}{||s_k||^2}\right\}+\mu \right)s_k.
\]
Then, the inequality $s_k^T \tilde{y}_k(\mu) > 0$ is always holds because
\[
  s_k^T \tilde{y}_k(\mu) = \max\{ 0,s_k^T y_k \}+\mu \|s_k\|^2 >0 .
\]

In the following, $\hat{H}_{k}(\mu)$ is the matrix constructed by the L-BFGS method using the initial matrix $
\hat{H}_k^{(0)}(\mu)$ and the vector pairs $(s_{k-i},\hat{y}_{k-i}(\mu)),~i=1,\cdots,m$, and the search direction 
is given as $d_k(\mu) = -\hat{H}_{k}(\mu)\nabla f(x_k)$.

The usual L-BFGS method uses $\gamma_k I$ as the initial matrix $H_k^{(0)}$, where $\gamma_k$ is a certain 
positive constant. Since $(B_k^{(0)})^{-1} = H_k^{(0)}$ and $\hat{H}_k^{(0)}$ is an approximation of 
$(B_k^{(0)} + \mu I)^{-1}$, we may set the initial matrix $\hat{H}_k^{(0)}(\mu)$ as
\begin{equation}
 \hat{H}_k^{(0)}(\mu) = (B_k^{(0)} + \mu I)^{-1} = \left(\frac{1}{\gamma_k} + \mu \right)^{-1} I = \frac{\gamma_k}
 {1+\gamma_k\mu}I.
\label{initial_hessian} 
\end{equation}

The proposed method generates the next iterate as $x_{k+1} = x_k + d_k(\mu)$ without a step length.
We control the parameter $\mu$ to guarantee the global convergence as in \cite{Ueda1}.
We exploit the idea of updating the trust-region radius in the TR-method to control $\mu$
to find an appropriate search direction, that is, we use the ratio of the reduction in the 
objective function value to that of the model function value.
We define a ratio function $r_k(d_k(\mu),\mu)$ by


\begin{eqnarray}
 r_k(d_k(\mu), \mu) = \frac{f(x_k) - f(x_k + d_k(\mu))}{f(x_k) - q_k(d_k(\mu),\mu)},
\end{eqnarray}
where $q_k : \mathbb{R}^n\times\mathbb{R} \rightarrow \mathbb{R}$ is given by
\[
  q_k(d_k(\mu),\mu) = f(x_k) + \nabla f(x_k)^Td_k(\mu) + \frac{1}{2}d_k(\mu)^T\hat{H}_k(\mu)^{-1}d_k(\mu).
\]
Note that we do not have to compute the matrix $\hat{H}_k(\mu)^{-1}$ explicitly in $q_k(d_k,
\mu)$.
Since $d_k(\mu) = -\hat{H}_k(\mu)\nabla f(x_k)$, we have $d_k(\mu)^T\hat{H}_k(\mu)^{-1}
d_k(\mu) = -d_k(\mu)^T\nabla f(x_k)$.
If the ratio $r_k(d_k(\mu), \mu)$ is large, i.e., the reduction in the objective function $f$ is 
sufficiently large compared to that of the model function,
we adopt $d_k(\mu)$ and decrease the parameter $\mu$. On the other hand, if $r_k(d_k,\mu)
$ is small, i.e., $f(x_k) - f(x_k+d_k)$ is small, we increase $\mu$ and compute $d_k(\mu)$ 
again.

Based on the above ideas, we propose the following regularized L-BFGS method.


\noindent\makebox[\linewidth]{\rule{\linewidth}{0.4pt}}


\begin{alg}
\label{alg_rlbfgs}
\centering{\textbf{\rm Regularized L-BFGS }}\\
\begin{description}
 \rm   \item[\textbf{Step 0}] Choose the parameters $\mu_0, \mu_{min},\gamma_1,
    \gamma_2,\eta_1,\eta_2,m$ such that $0<\mu_{min} \leq \mu_0, 0<\gamma_1
    \leq 1 <\gamma_2, 0<\eta_1<\eta_2\leq 1$ and $m > 0$. Choose initial
    point $x_0 \in \R^n$ and an initial matrix $\hat{H}^{0}_{k}$. Set 
    $k:= 0$.
    \item [\textbf{Step 1}] If some stopping criteria are satisfied, then 
    terminate. Otherwise go to step 2.
    \item[\textbf{Step 2}] \quad
            \begin{description}
         \item [\textbf{Step 2-0}] Set $l_k:=0$ and
         $\bar{\mu}_{l_k}=\mu_k.$
          \item[\textbf{Step 2-1}] 
          Compute $d_k(\bar{\mu_{l_k}})$ using Algorithm \ref{alg_twoloopreq}.
           \item[\textbf{Step 2-2}] Compute
           $r_k(d_k(\bar{\mu}_{l_k}),\bar{\mu}_{l_k})$. If
           $r_k(d_k(\bar{\mu}_{l_k}),\bar{\mu}_{l_k})< \eta_1,$ then update
           $\bar{\mu}_{l_{k+1}} = \gamma_2 \bar{\mu}_{l_{k}} $, set 
           $l_k = l_k +1, $ and  go to Step 2-1. Otherwise, go to Step $3$.
            \end{description}
    \item[\textbf{Step 3}] If $\eta_1\leq
    r_k(d_k(\bar{\mu}_{l_k}),\bar{\mu}_{l_k}) < \eta_2$ then update 
    $\mu_{k+1} = \bar{\mu}_{l_k}$.\\
    If  $r_k(d_k(\bar{\mu}_{l_k}),\bar{\mu}_{l_k}) \geq \eta_2$ then 
    update $\mu_{k+1} = \text{max}[\mu_{\text{min}}, \gamma_1
    \bar{\mu_{l_k}}]$. Update $x_{k+1} = x_k + d_k(\bar{\mu}_{l_k}).
     $ Set $k=k+1$ and go to Step 1.
\end{description}
\end{alg}

\noindent\makebox[\linewidth]{\rule{\linewidth}{0.4pt}}
\vspace*{0.15cm}
\newline In Step 2-1 we compute $\dk$ from $(s_k, \yk)$ by the \lb\ updating scheme in 
\cite{Nocedal}. The details of step 2-1 are  given as follows.
\vspace*{0.25cm}
\newline \noindent\makebox[\linewidth]{\rule{\linewidth}{0.4pt}}
\pagebreak
\begin{alg}
\label{alg_twoloopreq}
 \centering{\textbf{ \rm \lb\ with }$(s_k,\yk)$}
\begin{description}
\rm    \item[\textbf{Step 0}] Set $p  \leftarrow \gk.$
    
    \item[\textbf{Step 1}] Repeat the following process with 
    $i=k-1, k-2,\ldots, k-l;$
    \begin{align*}
        r_i & \leftarrow \tau_i s^T_i p,\\
        p & \leftarrow p - r_i(y_k +\mu s_k),\\
    \end{align*}
    where $\tau_i = (s^T_i(y_k +\mu s_k))^{-1}.$
    \item[\textbf{Step 2}] Set $q \leftarrow \hat{H}^{0}_k(\mu) p.$
    
    \item[\textbf{Step 3}] Repeat the following process with 
    $i=k-t,k-t+1,\ldots, k-1;$
    \begin{align*}
        \beta & \leftarrow \tau_i (y_k + \mu s_k)^T q, \\
        q & \leftarrow q + (r_i - \beta)s_i.\\
    \end{align*}
    \item[\textbf{Step 4}] Get the search direction by $\dk = -q. $
\end{description}
\end{alg}
 
\noindent\makebox[\linewidth]{\rule{\linewidth}{0.4pt}}
\vspace*{0.15cm}
\newline It is important to  note that  when $\mu$ varies, the regularized \lb\  does 
not have to store $\yk$ because the \lb stores $s_k$ and $y_k$ explicitly, and thus we 
can get $\yk$ immediately.

\section{Global convergence}

In this section, we show the global convergence of the proposed algorithm. To this end, we need the following 
assumptions.
\begin{asm}
 \label{assumption}
 \mbox{}
 \begin{description}
  \item {{\rm (i)}} The objective function $f$ is twice continuously differentiable.
  \item {{\rm (ii)}} The level set of $f$ at the initial point $x_0$ is compact, i.e., $\Omega=\{x\in\mathbb{R}^n|
  f(x)\leq f(x_0)\}$ is compact.
  \item {{\rm (iii)}} There exist positive constants $M_1$ and $M_2$ such that
      \[ M_1\|z\|^2 \leq z^T \nabla^2f(x) z \leq M_2\|z\|^2 ~\forall x \in \Omega \text{~ and ~} z \in \R^n.\]
  \item {{\rm (iv)}} There exists a minimum $f_{\min}$ of $f$.
  \item {{\rm (v)}} There exists a constant $\underline{\gamma}$ such that $\gamma_k \geq \underline{\gamma} > 0$ 
  for all $k$, where $\gamma_k$ is a parameter in \eqref{initial_hessian}.
 \end{description}
\end{asm}
The above assumptions are the same as those for the global convergence of the original L-BFGS method \cite{Liu}.

Under these assumptions, we have the following several properties.
First, let
\begin{eqnarray}
G(x) = \nabla^2 f(x),~~ G_k = G(x_k),~~
\bar{G}_k = \int_0^1 G(x_k+\tau s_k) d\tau . \nonumber
\end{eqnarray}
It then follows from  Taylor's theorem that
\begin{eqnarray}
f(x_k+d_k(\mu)) = f(x_k) + \nabla f(x_k)^T d_k(\mu) +\frac{1}{2} \int_0^1 d_k(\mu)^TG(x_k+\tau d_k(\mu))d_k(\mu) 
d\tau. \nonumber
\end{eqnarray}
Furthermore,
since $s_k = x_{k+1} - x_k $ and $y_k = \nabla f(x_{k+1}) - \nabla f(x_k) $, we have
\begin{eqnarray}
y_k=\bar{G}_ks_k, \label{yGs}
\end{eqnarray}
and hence we have
\begin{equation}
\hat{y}_k(\mu)= y_k + \mu s_k = (\bar{G}_k+\mu I)s_k.\label{yGms}
\end{equation}
It follows from {\rm Assumption} \ref{assumption} (iii) that $\lambda_{\min}(\bar{G}_k) \geq M_1$ and $
\lambda_{\max}(\bar{G}_k)\leq M_2$.
Therefore, we have that
\begin{eqnarray}
M_1\|s\|^2 &\leq& s_k^Ty_k \leq M_2\|s_k\|^2, \nonumber \\
\frac{1}{M_2}\|y_k\|^2 &\leq& s_k^Ty_k \leq \frac{1}{M_1}\|y_k\|^2, \label{jyokagen} \\
(M_1 + \mu)\|s_k\|^2 &\leq& s_k^T\hat{y}_k(\mu) \leq (M_2+\mu)\|s_k\|^2. \nonumber
\end{eqnarray}

Since the sequence $\{ x_k\}$ is included in the compact set $\Omega$ and $f$ is twice continuously differentiable 
under Assumption \ref{assumption} {\rm (i) and (ii)}, there exists a positive constant $L_f$ such that
\begin{eqnarray}
 \|\nabla f(x_k)\| \leq L_f~~~{\rm for~all}~k. \label{Lipschitz_f}
\end{eqnarray}

Now, we investigate the behavior of the eigenvalues of $\hat{B}_k(\mu)$, which is the inverse of $\hat{H}_k(\mu)$.
Note that the matrix $\hat{B}_k(\mu)$ is constructed by the BFGS formula
with vector pairs ($s_k,\hat{y}_k(\mu)$) and initial matrix $\hat{B}_k^{(0)}(\mu) = \hat{H}_k^{(0)}(\mu)^{-1}$. 
Thus, we have
\begin{equation}
\begin{split}
\hat{B}_k(\mu) &= \hat{B}_k^{(\tilde{m}_k)}(\mu) \\
\hat{B}_k^{(l+1)}(\mu) &= \hat{B}_k^{(l)}(\mu) -
\frac{\hat{B}_k^{(l)}(\mu) s_{j_l} s_{j_l}^T
\hat{B}_k^{(l)}(\mu)}{s_{j_l}^T \hat{B}_k^{(l)}(\mu) s_{j_l}} +
\frac{y_{j_l} y_{j_l}^T}{y_{j_l}^T s_{j_l}},~~l=0,\cdots,\tilde{m}_k-1
\end{split}
 \label{lbfgs2}
\end{equation}
where $\tilde{m}_k=\min\{k+1,m\}$ and $j_l = k-\tilde{m}+l$.
Note that these expressions are used in \cite{Byrd,Liu}.

We now focus on the trace and determinant of $\hat{B}_k(\mu)$.
First, we show that the trace of $\hat{B}_k^{(l)}(\mu)$ is $O(\mu)$.


\begin{lem}
\label{lem_tr}
Suppose that Assumption \ref{assumption} holds. Then,
\begin{eqnarray}
{\rm tr}(\hat{B}_{k}^{(l)}(\mu)) \leq M_3 + (2m+n)\mu,~~l=0,\cdots,~\tilde{m}_k \nonumber
\end{eqnarray}
where $  M_3 = \frac{n}{\underline{\gamma}}+m M_2.$

\end{lem}

\begin{prf}
{\rm
We have from {\rm Assumption \ref{assumption}}, {\rm (\ref{yGs}), and (\ref{jyokagen})} that
\begin{eqnarray}
\frac{\|\hat{y}_k(\mu)\|^2}{s_k^T\hat{y}_k(\mu)}
&=& \frac{ \|y_k\|^2 + 2\mu s_k^Ty_k + \mu^2\|s_k\| }{s_k^Ty_k + \mu \|s_k\|^2 } \nonumber \\
&=& \frac{\|y_k\|^2 + \mu s_k^Ty_k }{s_k^Ty_k + \mu \|s_k\|^2 } + \frac{\mu (s_k^Ty_k + \mu \|s_k\|^2)}{ s_k^Ty_k + 
\mu \|s_k\|^2} \nonumber \\
& \leq & \frac{\|y_k\|^2 + \mu s_k^Ty_k}{s_k^Ty_k} + \mu \nonumber \\
& \leq & \frac{\|y_k\|^2}{\frac{1}{M_2} \|y_k\|^2} + 2\mu \nonumber \\
&=& M_2 + 2\mu.  \label{yy_sy}
\end{eqnarray}
From the updating formula $\eqref{lbfgs2}$ of matrix $\hat{B}_k(\mu)$,
\begin{eqnarray}
{\rm tr}(\hat{B}_{k}^{(l)}(\mu))
={\rm tr}(\hat{B}_{k}^{(0)}(\mu)) + \sum_{t=0}^{l-1} \left(-
\frac{\| \hat{B}_k^{(t)}(\mu) s_{j_t} \|^2}{s_{j_t}^T \hat{B}_k^{(t)}(\mu) s_{j_t}}
+\frac{\|\hat{y_{j_t}}(\mu)\|^2}{s_{j_t}^T\hat{y}_{j_t}(\mu)} \right). \nonumber
\end{eqnarray}
It then follows from \eqref{yy_sy} that
\begin{eqnarray}
{\rm tr}(\hat{B}_{k}^{(l)}(\mu))
&\leq& {\rm tr}(\hat{B}_{k}^{(0)}(\mu)) + \sum_{t=0}^{l-1} \frac{\|\hat{y_{j_t}}(\mu)\|^2}{s_{j_t}^T\hat{y}_{j_t}(\mu)} \nonumber \\
&\leq& {\rm tr}(\hat{B}_{k}^{(0)}(\mu)) + l ( M_2+ 2\mu ) \nonumber \\
&=& n\left(\frac{1}{\gamma_{k}}+\mu\right) + l (M_2 + 2\mu) \nonumber \\
& \leq & n\left(\frac{1}{\underline{\gamma}}+\mu\right) + m (M_2 + 2\mu) \nonumber \\
& \leq &  M_3 +(2m+n)\mu. \nonumber
\end{eqnarray}
This completes the proof.}
\qed
\end{prf}


The next lemma gives a lower bound for the determinant of $\hat{B}_k(\mu)$.
\begin{lem}
\label{lem_det}
Suppose that Assumption \ref{assumption} holds. Then,
\begin{eqnarray}
{\rm det}(\hat{B}_{k}(\mu) ) \geq M_4\mu^n,\nonumber
\end{eqnarray}
where
\[
  M_4 = \left(\frac{1}{2m+n}\right)^{m}.
\]
\end{lem}
\begin{prf}
{\rm
Note that the determinant of the approximate matrix updated by the BFGS updating scheme has
the following property {\rm \cite{Pearson,Powell}}:
\begin{eqnarray}
{\rm det}(\hat{B}_k^{(l+1)}(\mu))
= {\rm det}(\hat{B}_k^{(l)}(\mu))\frac{s_{j_l}^T \hat{y}_{j_l}(\mu)}{s_{j_l}^T \hat{B}_{k-1}^{(l)}(\mu) s_{j_l} }. \nonumber
\end{eqnarray}
Then, we have
\begin{eqnarray}
{\rm det}(\hat{B}_{k}(\mu)) &=& {\rm det}(\hat{B}_{k}^{(0)}(\mu))
\prod_{l=0}^{\tilde{m}-1}\frac{ s_{j_l}^T\hat{y}_{j_l}(\mu) }{ s_{j_l}^T\hat{B}_{k-1}^{(0)}(\mu) s_{j_l} } \nonumber \\
&=&{\rm det}( \hat{B}_{k}^{(0)}(\mu) ) \prod_{l=0}^{\tilde{m}-1}\frac{ s_{j_l}^T\hat{y}_{j_l}(\mu) }{ s_{j_l}^Ts_{j_l} }\frac{ s_{j_l}^Ts_{j_l} }{ s_{j_l}^T\hat{B}_{k}^{(l)}(\mu)s_{j_l} } \nonumber \\
&\geq&{\rm det}( \hat{B}_{k}^{(0)}(\mu)) \prod_{l=0}^{\tilde{m}-1}\frac{ s_{j_l}^T\hat{y}_{j_l}(\mu) }{ s_{j_l}^Ts_{j_l} }\frac{ s_{j_l}^Ts_{j_l} }{ \lambda_{\max}(\hat{B}_{k}^{(l)}(\mu)) s_{j_l}^Ts_{j_l} } \nonumber \\
&=& {\rm det}(\hat{B}_k^{(0)}(\mu)) \prod_{l=0}^{\tilde{m}-1}\frac{s_{j_l}^T\hat{y}_{j_l}(\mu)}{\|s_{j_l}\|^2}\frac{1}{\lambda_{\max}(\hat{B}_k^{(l)}(\mu))}. \nonumber
\end{eqnarray}
Since $B_k^{(0)}(\mu)$ is symmetric positive-definite,
{\rm Lemma} $\ref{lem_tr}$ implies that $\lambda_{\max}(\hat{B}_{k}^{(l)}(\mu)) \geq M_3+(2m+n)\mu$.
Furthermore, we have $\frac{s_{j_l}^T\hat{y}_{j_l}(\mu)}{\|s_{j_l}\|^2} \geq M_1+\mu$ from $\eqref{jyokagen}$. Therefore, it follows that
\begin{eqnarray}
{\rm det}(\hat{B}_k(\mu)) &\geq& {\rm det}( \hat{B}_{k}(\mu)^{(0)} )
\left( \frac{M_1+\mu}{M_3+(2m+n)\mu} \right)^{\tilde{m}} \nonumber \\
&=&{\rm det}\left( \frac{1+\gamma_{k}\mu}{\gamma_{k}}I\right)\left( \frac{M_1+\mu}{M_3+(2m+n)\mu}\right)^{\tilde{m}} \nonumber \\
& \geq & \left( \frac{1}{\gamma_{k}}+\mu \right)^n \left(\frac{1}{2m+n}\right)^{\tilde{m}} \nonumber \\
& \geq & \left( \frac{1}{2m+n}\right)^m \mu^n \nonumber \\
&=& M_4 \mu^n. \nonumber
\end{eqnarray}
This completes the proof.}
\qed
\end{prf}

From the above two lemmas, we have $\lambda_{\max}(\hat{H}_{k}(\mu) ) \rightarrow 0$ as $
\mu \rightarrow \infty$.


\begin{lem}
\label{lem_H0}
Suppose that Assumption \ref{assumption} holds. Then, for all $k\geq 0$,
\begin{eqnarray}
\lambda_{\max}(\hat{H}_k(\mu)) \leq M_5 \frac{1}{\mu},~\forall \mu \in [ \mu_{\min},\infty),\nonumber
\end{eqnarray}
where
\[
M_5 =\frac{1}{M_4 n^{n-1}}\left( \frac{M_3}{\mu_{\min}} + (2m+n) \right)^{n-1}.
\]
Furthermore,
$\lim_{\mu \to \infty} \lambda_{\max}(\hat{H}_k(\mu)) = 0.$
\end{lem}
\begin{prf}
{\rm
We have from {\rm Lemmas} $\ref{lem_tr}$ and $\ref{lem_det}$ that
\begin{eqnarray}
{\rm tr}(\hat{B}_k(\mu)) &\leq& M_3 + (2m+n)\mu, \nonumber \\
{\rm det}(\hat{B}_{k}(\mu) ) &\geq& M_4\mu^n. \nonumber
\end{eqnarray}
Since $\hat{B}_k(\mu)$ is symmetric positive-definite, we have
\begin{eqnarray}
{\rm tr}(\hat{B}_{k}(\mu)) &\geq& \lambda_{\min}(\hat{B}_k(\mu)) \nonumber \\
 \det(\hat{B}_{k}(\mu)) &\leq&  \lambda_{\min}(\hat{B}_k(\mu) )\{\lambda_{\max}(\hat{B}_k(\mu))\}^{n-1}. \nonumber
\end{eqnarray}
Therefore, we have
\begin{eqnarray}
\lambda_{\min}(\hat{B}_k(\mu)) &\geq& \frac{ \det(\hat{B}_{k}(\mu)) }{ \{\lambda_{\max}
(\hat{B}_k(\mu))\}^{n-1} } \nonumber \\
 & \geq & \frac{ M_4 \mu^n }{ \{ (M_3 +(2m+n)\mu)\}^{n-1} }. \nonumber
\end{eqnarray}
It then follows from Assumption {\rm \ref{assumption} (v)} that
\begin{eqnarray}
\lambda_{\max}(\hat{H}_k(\mu)) &=& \frac{1}{\lambda_{\min}(\hat{H}_k^{-1}(\mu))} \nonumber \\
&=&\frac{1}{\lambda_{\min}(\hat{B}_k(\mu))} \nonumber \\
& \leq & \frac{ \{(M_3 +(2m+n)\mu) \}^{n-1} }{ M_4\mu^n } \nonumber \\
&=& \frac{1}{M_4 }\left( \frac{M_3 + (2m+n)\mu }{\mu}\right)^{n-1}\frac{1}{\mu}. \label{lambda_kagen}
\end{eqnarray}
Since $\mu\geq \mu_{\min}$, we have
\[
 \frac{M_3 + (2m+n)\mu}{\mu} = \frac{M_3}{\mu}+(2m+n) \leq \frac{M_3}{\mu_{\min}}+(2m+n).
\]
It then follows from \eqref{lambda_kagen} that
\begin{eqnarray}
\lambda_{\max}(\hat{H}_k(\mu))
& \leq & \frac{1}{M_4}\left( \frac{M_3}{\mu_{\min}} + (2m+n) \right)^{n-1} \frac{1}{\mu} \nonumber \\
& = & M_5 \frac{1}{\mu}. \nonumber
\end{eqnarray}
Hence, we have
\[\lim_{\mu \to \infty} \lambda_{\max}(\hat{H}_k(\mu)) = 0.  \]
This completes the proof.}
\qed
\end{prf}

Now, we give an upper bound for $\|d_k(\mu)\|$.


\begin{lem}
\label{lem_upd}
Suppose that Assumption \ref{assumption} holds. Then,
\begin{eqnarray}
\|d_k(\mu)\| \leq U_d, \nonumber
\end{eqnarray}
where
\[
 U_d =  \frac{L_f M_5}{\mu_{\min}}.
\]
\end{lem}
\begin{prf}
{\rm
From the definition of $d_k(\mu)$, \eqref{Lipschitz_f}, and Lemma \ref{lem_H0}, we have that
\begin{eqnarray}
\|d_k(\mu)\|&=&\|\hat{H}_k(\mu)\nabla f(x_k)\| \nonumber \\
&\leq& \|\hat{H}_k(\mu)\| \|\nabla f(x_k)\| \nonumber\\
&=&\lambda_{\max}(\hat{H}_k(\mu)) \|\nabla f(x_k)\| \nonumber \\
& \leq & \lambda_{\max}(\hat{H}_k(\mu)) L_f \nonumber \\
& \leq & \frac{L_f M_5}{\mu} \nonumber \\
& \leq & \frac{L_f M_5}{\mu_{\min}} = U_d. \nonumber
\end{eqnarray}
This completes the proof.}
\qed
\end{prf}

Lemma \ref{lem_upd} implies that
\[
 x_k + \nu d_k(\mu) \in \Omega + {\rm B}(0,U_d),~~~\forall \nu \in [0,1],~~~\forall \mu \in [\mu_{\min},\infty),~~~\forall k \geq 0.
\]
Moreover, since $\Omega+{\rm B}(0,U_d)$ is compact and $f$ is twice continuously differentiable, $\nabla f(x_k)$ is Lipschitz continuous on
$\Omega + {\rm B}(0,U_d)$. That is, there exists a positive constant $L_g$ such that
\begin{equation}
 \|\nabla^2 f(x_k)\| \leq L_g~~~\forall x_k \in \Omega+{\rm B}(0,U_d). \label{Lipshitz_g}
\end{equation}

Next, we investigate the values of $\mu$ that satisfy the termination condition 
$r_k(d_k(\mu),$ $\mu)\geq \eta_1$ in the inner iterations of Step 2-2 in Algorithm
 \ref{alg_rlbfgs}.


\begin{lem}
\label{lem_ratio}
Suppose that Assumption \ref{assumption} holds. Then, we have
\begin{multline*}
f(x_k) - f(x_k+d_k(\mu)) - \eta_1 (f(x_k) - q_k(d_k(\mu),\mu)) \\
 \geq \frac{1}{2}((2-\eta_1) \lambda_{\min}(\hat{H}_k(\mu)^{-1}) - L_g)\|d_k(\mu)\|^2.
\end{multline*}
\end{lem}
\begin{prf}
{\rm
We have from  Taylor's theorem that
\[ \begin{aligned}[t]
f(x_k+d_k(\mu)) & =  f(x_k) + \int_0^1 \nabla f(x_k+\tau d_k(\mu))^Td_k(\mu)d\tau 
\nonumber \\
   & =  f(x_k) + \nabla f(x_k)^Td_k(\mu) + \int_0^1(\nabla f(x_k+\tau d_k(\mu)) \nonumber
   \\
&  \qquad \qquad \qquad\qquad  \qquad \qquad \qquad  \quad -\nabla f(x_k))^Td_k(\mu) d
\tau.\nonumber
\end{aligned} \]
From the Lipschitz continuity of $\nabla f(x_k)$ in (\ref{Lipshitz_g}), we get

\begin{multline*}
f(x_k) - f(x_k+d_k(\mu)) - \eta_1 (f(x_k) - q_k(d_k(\mu),\mu)) \nonumber\\
=-\nabla f(x_k)^Td_k(\mu) - \int_0^1(\nabla f(x_k+\tau d_k(\mu))-
\nabla f(x_k))^Td_k(\mu) d\tau \nonumber\\
- \frac{\eta_1}{2}d_k(\mu)^T(\hat{H}_k(\mu)^{-1})d_k(\mu) \nonumber
\end{multline*}
\[ \begin{aligned}[t]
&=\frac{(2-\eta_1)}{2}d_k(\mu)^T(\hat{H}_k(\mu)^{-1})d_k(\mu) - \int_0^1(\nabla f(x_k+
\tau d_k(\mu)) \nonumber\\
&   \qquad \qquad\qquad  \qquad \qquad \qquad \qquad \qquad \qquad \qquad \qquad
\qquad -\nabla f(x_k))^Td_k(\mu) d\tau \nonumber\\
&\geq\frac{(2-\eta_1)}{2}\lambda_{\min}(\hat{H}_k(\mu)^{-1})\|d_k(\mu)\|^2-\int_0^1L_g\tau \|d_k(\mu)\|^2d\tau \nonumber\\
&=\frac{1}{2}((2-\eta_1)\lambda_{\min}(\hat{H}_k(\mu)^{-1}) - L_g)\|d_k(\mu)\|^2. \nonumber
\end{aligned} \]
This completes the proof.}
\qed
\end{prf}

From Lemma \ref{lem_ratio},
if $\mu$ satisfies
\begin{eqnarray}
\lambda_{\min}( \hat{H}_k^{-1}(\mu) )\geq \frac{L_g}{2-\eta_1}, \label{step2fin}
\end{eqnarray}
then we have
\begin{equation}
 r_k(d_k(\mu),\mu)\geq \eta_1, \label{lem_rkagen}
\end{equation}
that is, the inner loops of Algorithm \ref{alg_rlbfgs} must terminate.

Next, we give an upper bound for the parameter $\mu_k$.


\begin{lem}
\label{lem_mu}
Suppose that Assumption \ref{assumption} holds. Then, for any $k \geq 0$,
\begin{eqnarray}
\mu_k^{\ast}\leq U_{\mu}, \nonumber
\end{eqnarray}
where
\[
 U_{\mu} =\gamma_2 M_5 \frac{L_g}{2-\eta_1}.
\]
\end{lem}
\begin{prf}
{\rm
If $\bar{\mu}_{l_k}$ satisfies \eqref{step2fin}, then $r_k(d_k(\bar{\mu}_{l_k}),\bar{\mu}_{l_k})\geq\eta_1$ from {\rm Lemma} $\ref{lem_ratio}$.
Therefore, the inner loops must terminate, and we set $\mu_k^{\ast} =\bar{\mu}_{l_k}$.

Now, we give the termination condition on $\mu$ for the inner loop.
We have from {\rm Lemma {\ref{lem_H0}}} that
\begin{eqnarray}
 \lambda_{\min}(\hat{H}_k^{-1}(\mu))
 &=& \frac{1}{\lambda_{\max}(\hat{H}_k(\mu))} \nonumber \\
 &\geq& \frac{\mu}{M_5}. \label{muM5}
\end{eqnarray}
It then follows from {\rm (\ref{step2fin})} that
the termination condition of the inner loop holds when
\begin{eqnarray}
 \mu \geq M_5\frac{L_g}{2-\eta_1}. \label{step2fin2}
\end{eqnarray}
Note that if the inner loop terminates at ${l_k}$, then \eqref{step2fin2} does not hold
with $\mu = \mu_{l_k-1}$, that is,
\begin{eqnarray}
 \bar{\mu}_{l_k-1} < M_5 \frac{L_g}{2-\eta_1}.\nonumber
\end{eqnarray}
Since $\mu_k^{\ast} = \bar{\mu}_{l_k} = \gamma_2 \bar{\mu}_{l_k-1}$, we have
\begin{eqnarray}
\mu_k^{\ast} = \gamma_2 \bar{\mu}_{l_k-1} < \gamma_2 M_5 \frac{L_g}{2-\eta_1} = U_{\mu}.
\end{eqnarray}
This completes the proof.
}
\qed
\end{prf}

Next, we give a lower bound for the reduction in the model function $q_k$.


\begin{lem}
\label{lem_modelkagen}
Suppose that Assumption \ref{assumption} holds. Then, we have
\begin{eqnarray}
f(x_k)-q_k(d_k(\mu),\mu)&\geq& M_6 \|\nabla f(x_k)\|^2, \nonumber
\end{eqnarray}
where
\[
 M_6 = \frac{1}{2(M_3+(2m+n)\mu_{\min})}.
\]
\end{lem}
\begin{prf}
{\rm
It follows from the definition of the model function $q_k(d_k(\mu),\mu)$ and {\rm Lemmas
\ref{lem_tr}, \ref{lem_mu}} that
\begin{eqnarray}
f(x_k)-q_k(d_k(\mu),\mu) &=& -\frac{1}{2} d_k(\mu)^T ( \hat{H}_k^{-1}(\mu) )d_k(\mu) - \nabla f(x_k)^Td_k(\mu) 
\nonumber \\
&=& -\frac{1}{2}\nabla f(x_k)^T\hat{H}_k(\mu) \nabla f(x_k) + \nabla f(x_k)^T \hat{H}_k(\mu)\nabla f(x_k) \nonumber 
\\
&=& \frac{1}{2} \nabla f(x_k)^T\hat{H}_k(\mu) \nabla f(x_k)\nonumber \\
&\geq& \frac{1}{2} \lambda_{\min}( \hat{H}_k(\mu) )\|\nabla f(x_k)\|^2 \nonumber \\
&=& \frac{\|\nabla f(x_k)\|^2}{2\lambda_{\max}( \hat{H}_k^{-1}(\mu) )} \nonumber \\
&=& \frac{\|\nabla f(x_k)\|^2}{2\lambda_{\max}( \hat{B}_k(\mu) )} \nonumber \\
&\geq& \frac{\|\nabla f(x_k)\|^2}{2{\rm tr}( \hat{B}_k(\mu) )} \nonumber \\
&\geq& \frac{\|\nabla f(x_k)\|^2}{2(M_3+(2m+n)\mu)} \nonumber\\
&=& \frac{1}{2(M_3+(2m+n)U_{\mu})} \|\nabla f(x_k)\|^2 \nonumber \\
&=& M_6 \|\nabla f(x_k)\|^2.\nonumber
\end{eqnarray}
This completes the proof.
}
\qed
\end{prf}

From this lemma, we can give a lower bound for the reduction in the objective function value when $x_k$ is not a 
stationary point.


\begin{lem}
\label{lem_fkagen}
Suppose that Assumption \ref{assumption} holds.
If there exists a positive constant $\epsilon_g$ such that $\|\nabla f(x_k)\| \geq \epsilon_g$, then we have 
$f(x_k)-f(x_{k+1}) \geq \rho\epsilon_g^2$, where $\rho = \eta_1 M_6$.
\end{lem}
\begin{prf}
{\rm
It follows from {\rm Lemmas \ref{lem_mu} and \ref{lem_modelkagen}} that
\begin{eqnarray}
f(x_k)-f(x_{k+1}) &\geq& \eta_1(f(x_k)-q_k(d_k(\mu_k^{\ast}),\mu_k^{\ast})) \nonumber \\
&\geq& \eta_1 M_6 \|\nabla f(x_k)\|^2 \nonumber \\
&\geq&\rho\epsilon_g^2. \nonumber
\end{eqnarray}
This completes the proof.
}
\qed
\end{prf}

We are now in a position to prove the main theorem of this section.


\begin{thm}
Suppose that Assumption \ref{assumption} holds. Then,
$\liminf_{k\to\infty}\|\nabla f(x_k) \|=0$
 or there exists $K\geq0$ such that $\|\nabla f(x_K)\|=0$.
\end{thm}
\begin{prf}
{\rm
Suppose the contrary, that is, there exists a positive constant $\epsilon_g$ such that
$\|\nabla f(x_k)\| \geq \epsilon_g$ for all $k \geq 0$.
It follows from {\rm Lemma} \ref{lem_fkagen} that
\begin{eqnarray}
f(x_0)-f(x_k) &=& \sum_{j=0}^{k-1} (f(x_j)-f(x_{j+1})) \nonumber \\
&\geq&\sum_{j=0}^{k-1} \rho \epsilon_g^2 \nonumber \\
&=& \rho \epsilon_g^2 k. \nonumber
\end{eqnarray}
Taking $k \to \infty$, the right-hand side of the final inequality goes to infinity, and hence
\[\lim_{k\to\infty} f(x_k) = -\infty.\]
This contradicts the existence of $f_{\min}$ in Assumption {\ref{assumption}} (iv).
This completes the proof.}
\qed
\end{prf}

\section{Implementation issues}


The regularized \lb\ method does not  use a line search, and hence it can
not take a longer step.
Moreover in our experience, the trust-region ratio  of the regularized \lb\ does not improve 
well and the regularized parameter $\mu$ becomes very large for some large-scale test 
problems.
Both cases result in a short step, and hence the method conducts a large number of iteration 
to reach a solution. To overcome this difficulty we propose two techniques in this section. We 
also discuss how to set $\gamma_k$ in  the initial matrix $\hk$.
\subsection{Simultaneous use with Wolfe line search}

The next iterate  with line search is given as
\begin{equation}
\label{iter}
 x_{k+1} = x_k + \alpha_k d_k,
\end{equation}
where $\alpha_k$ is a step length. The usual L-BFGS \cite{Liu,Nocedal} uses a step length $\alpha_k$ 
that satisfies the Wolfe conditions,
\begin{align}
&  f(x_k+\alpha_k d_k)   \leq f(x_k) + c_1 \alpha_k d_k^T\gk;\label{Armijo}\\
&  d_k^T\g f(x_k+ \alpha_k d_k)   \geq c_2 d_k^T \gk;\label{Curvature}\\
& |d_k^T \g f(x_k + \alpha_k d_k) |  \leq  c_2 |d_k^T \gk | \label{strongWOLFE},
\end{align}
where $0<c_1<c_2<1$. 
Note that $\alpha_k$ can be larger than 1. Thus, $\alpha_k d_k$ might be larger, and make a large reduction of 
$f$.
Thus, it might be reasonable to use a line search as well as the regularization technique. However, finding
$\alpha_k$ takes much time, and hence we must avoid it if $\alpha_k d_k$ does not enough 
improvement.








For the efficient use of the line search, we exploit 
curvature condition \eqref{Curvature}. It is known that the curvature condition 
ensures that the step is not too short. Therefore, after step 3 of Algorithm \ref{alg_rlbfgs}, we 
first check weather $x_{k}+\dk$ satisfies curvature condition ($\ak=1$ in \eqref{Curvature}) or 
not. The dissatisfaction of the curvature condition implies that $x_k +\dk $ is a short step. 
Thus, we  compute $\ak$ by the strong Wolfe condition so that we can take a longer step.
 More precisely, we search $\alpha_k$ from $x_k +\dk$ with the direction $\dk$ so that 
  \eqref{Armijo}-\eqref{strongWOLFE} hold with $x_k:=x_k+\dk$, $d_k:=\dk,$ and then set 
  $x_{k+1}=x_k+(1+\alpha_k) \dk$. 





We now discuss the conditions under which we conduct the Wolfe line search. As mentioned 
above, we exploit the strong Wolfe condition \eqref{Curvature} when the following 
conditions hold,
\begin{equation}
\dk^T\g f(x_k+ \dk) < c_2 \dk^T \gk \quad \text{ and }\quad \mu=\mu_{min}.
\label{RL-BFGS-SW_condition}
\end{equation}
Note that  $\| d_k(\mu_{min})\| \geq \| \dk \|$ if $\mu> \mu_{min} $. Thus, $d_k(\mu_{min})$ is
the largest step when we apply RL-BFGS only. Condition \eqref{RL-BFGS-SW_condition} implies whenever 
$x_k+ d_k(\mu_{min})$ to  make better progress, we take a longer step via a strong Wolfe line search.
We call this method regularized L-BFGS with strong 
Wolfe line search method (RL-BFGS-SW) as an extended version of the proposed method. Now, we 
propose the RL-BFGS-SW as follows.\\
\noindent\makebox[\linewidth]{\rule{\linewidth}{0.4pt}}\\
\begin{alg}
\centering{\textbf{\rm R\lb\  with line search (RL-BFGS-SW)}}\\
\begin{description}
\label{alg_rlbfgs_withlineser}
\rm  \item[\textbf{Step 0}] Choose the parameters $\mu_0, \mu_{min},\gamma_1,
    \gamma_2,\eta_1,\eta_2,m$ such that $0<\mu_{min} \leq \mu_0, 0<\gamma_1
    \leq 1 <\gamma_2, 0<\eta_1<\eta_2\leq 1$ and $m > 0$. Choose initial
    point $x_0 \in \R^n$ and an initial matrix $\hat{H}^{0}_{k}$. Set 
    $k:= 0$.
    \item[\textbf{Step 1}] If some stopping criteria are satisfied, then 
    terminate. Otherwise go to step 2.
\item[\textbf{Step 2}] \quad
\begin{description}
\item[\textbf{Step 2-0}] Set $l_k:=0$ and
         $\bar{\mu}_{l_k}=\mu_k$.
 \item[\textbf{Step 2-1}] 
          Compute $d_k(\bar{\mu_{l_k}})$ by Algorithm \ref{alg_twoloopreq}.
 \item[\textbf{Step 2-2}] Compute
           $r_k(d_k(\bar{\mu}_{l_k}),\bar{\mu}_{l_k})$. If
           $r_k(d_k(\bar{\mu}_{l_k}),\bar{\mu}_{l_k})< \eta_1,$ then update\\
           $\bar{\mu}_{l_{k+1}} = \gamma_2 \bar{\mu}_{l_{k}} $, set 
           $l_k = l_k +1, $ and  go to Step 2-1. Otherwise, go to Step $3$.            
\end{description}
\item[\textbf{Step 3}] If $\eta_1\leq r_k(d_k(\bar{\mu}_{l_k}),\bar{\mu}_{l_k})
      < \eta_2$ then update $\mu_{k+1} = \bar{\mu}_{l_k}$.\\
         If  $r_k(d_k(\bar{\mu}_{l_k}),\bar{\mu}_{l_k}) \geq \eta_2$ then 
    update $\mu_{k+1} = \text{max}[\mu_{\text{min}}, \gamma_1
    \bar{\mu_{l_k}}]$. 
    \item[\textbf{Step 4}] If $\dk^T\g f(x_k+ \dk) < c_2 \dk^T \gk$ 
    and $\mu_k=\mu_{min}$,\\
    then find $\alpha_k$ by 
    strong Wolfe line search and set $x_{k+1}  = x_k + d_k+ \alpha_k d_k$.\\
    Otherwise $x_{k+1}  = x_k + d_k$.\\
    Set $k=k+1$ and go to Step 1.
  
\end{description}
\end{alg}

\noindent\makebox[\linewidth]{\rule{\linewidth}{0.4pt}}
\vspace*{0.15cm}

Under this procedure, we must replace $s_k$ and $\yk$ whenever we use the strong Wolfe line 
search. We summarize  $y_k, \yk$ and $s_k$ in Table \ref{Table1}. 

\begin{table}[H]
\caption{Comparison of $y_k, \hat{y}_k$ and $s_k$.}
\centering
\begin{tabular}{|p{0.3\linewidth} |p{0.25\linewidth}| p{0.33\linewidth}| }
\hline
L-BFGS & RL-BFGS & RL-BFGS-SW \\
\ & \ & (when line search is used) \\
\hline
$x_{k+1}=x_k+\alpha_k d_k$     & $x_{k+1}=x_k+ \dk$      &  $x_{k+1}=x_k+(\alpha_k
+1) \dk$ \\
$s_k=\alpha_k d_k$                    &$s_k=\dk$                       & $s_k=(\alpha_k+1)\dk$  \\
$y_k=\g f(x_{k+1}) - \g f(x_{k})$ & $\yk=y_k+\mu_k \dk$   & $\yk=y_k+\mu_k (\alpha_k
+1)\dk $  \\
\hline 
\end{tabular}
 \label{Table1}
\end{table}

Note that we do not need to evaluate a new function and gradient values for conditions 
\eqref{RL-BFGS-SW_condition} since we have the gradient and function values to calculate the
 ratio $r_k(\dk, \mu_k)$. Note also that Algorithm \ref{alg_rlbfgs_withlineser} still has the 
 global convergence property since $f(x_k+\dk) > f(x_k + (1+\alpha_k)\dk)$ from \eqref{Armijo} and 
 \eqref{RL-BFGS-SW_condition}.

\vspace{5pt}
\noindent
\subsection{Nonmonotone decreasing technique}
In Algorithm \ref{alg_rlbfgs}, we control the regularized parameter $\mu$ to satisfy the descent condition 
$f(x_{k+1}) < f(x_k)$.
However, $\mu$ sometimes becomes quite large for some ill-posed problems.
In this situation, we require a large number of function evaluations.
Therefore, we use the concept of a nonmonotone line search technique \cite{Grippo,Sun} to overcome 
the difficulty. We replace the ratio function $r_k(d_k(\mu),\mu)$ with the following new ratio 
function $\bar{r}_k(d_k(\mu),\mu)$:
\begin{eqnarray}
 \bar{r}_k(d_k(\mu),\mu) = \frac{ \max_{0\leq j\leq m(k)}f(x_{k-j}) - f(x_k+d_k(\mu)) }{ f(x_k) - 
 q_k(d_k(\mu),\mu) }, \nonumber
\end{eqnarray}
where
\[
 m(0) = 0,~~0 \leq m(k) \leq \min\{m(k-1)+1,M\},
\]
and $M$ is a nonnegative integer constant. This modification retains the global convergence of 
the regularized L-BFGS method.

In the numerical experiments reported in the next section, when $k<M$, we use the original 
ratio function $r_k(d_k(\mu),\mu)$, and if $k\geq M$ then we use the new ratio function 
$\bar{r}_k(d_k(\mu),\mu)$.

\vspace{5pt}
\noindent
\subsection{Scaling initial matrix}
The regularized L-BFGS method uses the following initial matrix in each iteration:
\[
 \hat{H}_k^{(0)}(\mu) = \frac{\gamma_k}{1+\gamma_k\mu} I.
\]
The parameter $\gamma_k$ represents the scale of $\nabla^2f(x)$.
Thus, we exploit 
the scaling parameter $\gamma_k$ used in \cite{Barazilai,Burke2,Nocedal,Raydan,Shanno},
that is, we set
\[
 \gamma_k = \frac{s_{k-1}^Ty_{k-1}}{\|y_{k-1}\|^2}.
\]
It is known that the L-BFGS method with this scaling in the initial matrix has an efficient 
performance \cite{Burke2,Nocedal}.
Note that we require $\gamma_k>0$ to ensure the positive-definiteness of $\hat{H}_k^{(0)}(\mu)$.
If $s_{k-1}^Ty_{k-1} < \alpha \|s_{k-1}\|^2$, then we set $\gamma_k = \alpha\frac{\|s_{k-1}
\|^2}{\|y_{k-1}\|^2}$, where $\alpha$ is a small positive constant.


\section{Numerical results}
In this section, we compare the L-BFGS, the regularized L-BFGS (RL-BFGS), and the regularized L-BFGS
with line search (RL-BFGS-SW). For the regularized ones, we adopt the nonmonotone techniques
and the initial matrix discussed in Section 4. We have used 
MCSRCH (Line search routine) and parameters of the original L-BFGS~\cite{NocedalURL} to find a step length in 
the RL-BFGS-SW.


We have solved 313 problems chosen from CUTEst~\cite{Gould}. All algorithms
were coded in MATLAB 2018a. We have used Intel Core i5 1.8 GHz CPU with 8 GB RAM on Mac 
OS. We have chosen an initial point $x_0$ given in CUTEst. 

We set the same termination criteria as in the original L-BFGS, that is,
\begin{equation}
\frac{\|\gk\|}{\text{max}(1,\|x_k\|)} < 10^{-5}\  \text{or}\ n_f > 10000,
\end{equation}
where $n_f$ is the number of function evaluations. These criteria are similar to those in 
\cite{NocedalURL}.  We regard the trails as fail when $n_f > 10000$.


We compare the algorithms from the distribution function proposed in \cite{Dolan}. Let $
\mathcal{S}$ 
be a  set of solvers and let $\mathcal{P_S}$ be a set of problems that can be solved by all 
algorithms in $\mathcal{S}$. We measure required evaluations to solve problem $p$ by solver
$s\in \mathcal{S}$ as $t_{p,s}$, and the best $t_{p,s}$ for each $p$ as $t^*_p,$ which means
$t^*_p = \text{min}\{t_{p,s}| s \in \mathcal{S} \}$. The distribution function 
$F_s^{\mathcal{S}}(\tau)$, for a method $s$ is defined by,
\begin{equation}
 F_s^{\mathcal{S}}(\tau) = \frac{|\{p \in \mathcal{P}_{\mathcal{S}}| t_{p,s} \leq \tau t^*_p
 \}|}{|\mathcal{P}_{\mathcal{S}}|}, \ \tau \geq 1.
\end{equation}
The algorithm whose $F_s^{\mathcal{S}}(\tau)$ is close to 1 is considered  to be 
superior compare to other algorithm in $\mathcal{S}$.


\subsection{Numerical behavior for some parameters in R\lb}

Since the R\lb uses several parameters, we need to investigate the effect of 
these parameters so that we choose optimal ones.

\begin{figure}[H]
\centering
\begin{minipage}[b]{.5\textwidth}
\centering
\includegraphics[scale=0.12]{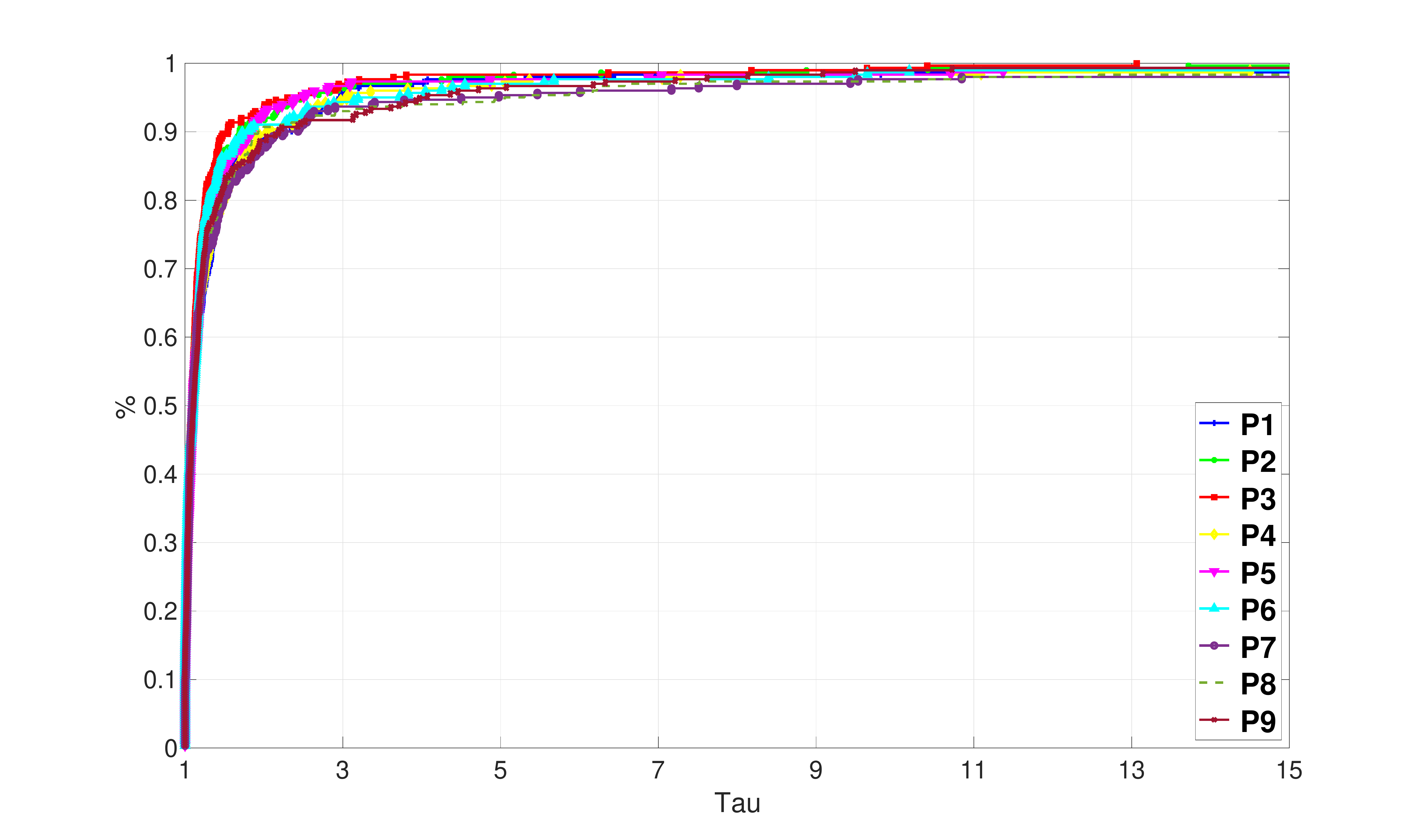}
 \caption{Comparison of  $(\gamma_1,\gamma_2)$.}
 \label{figgama}
\end{minipage}%
\begin{minipage}[b]{.5\textwidth}
\centering
\includegraphics[scale=.12]{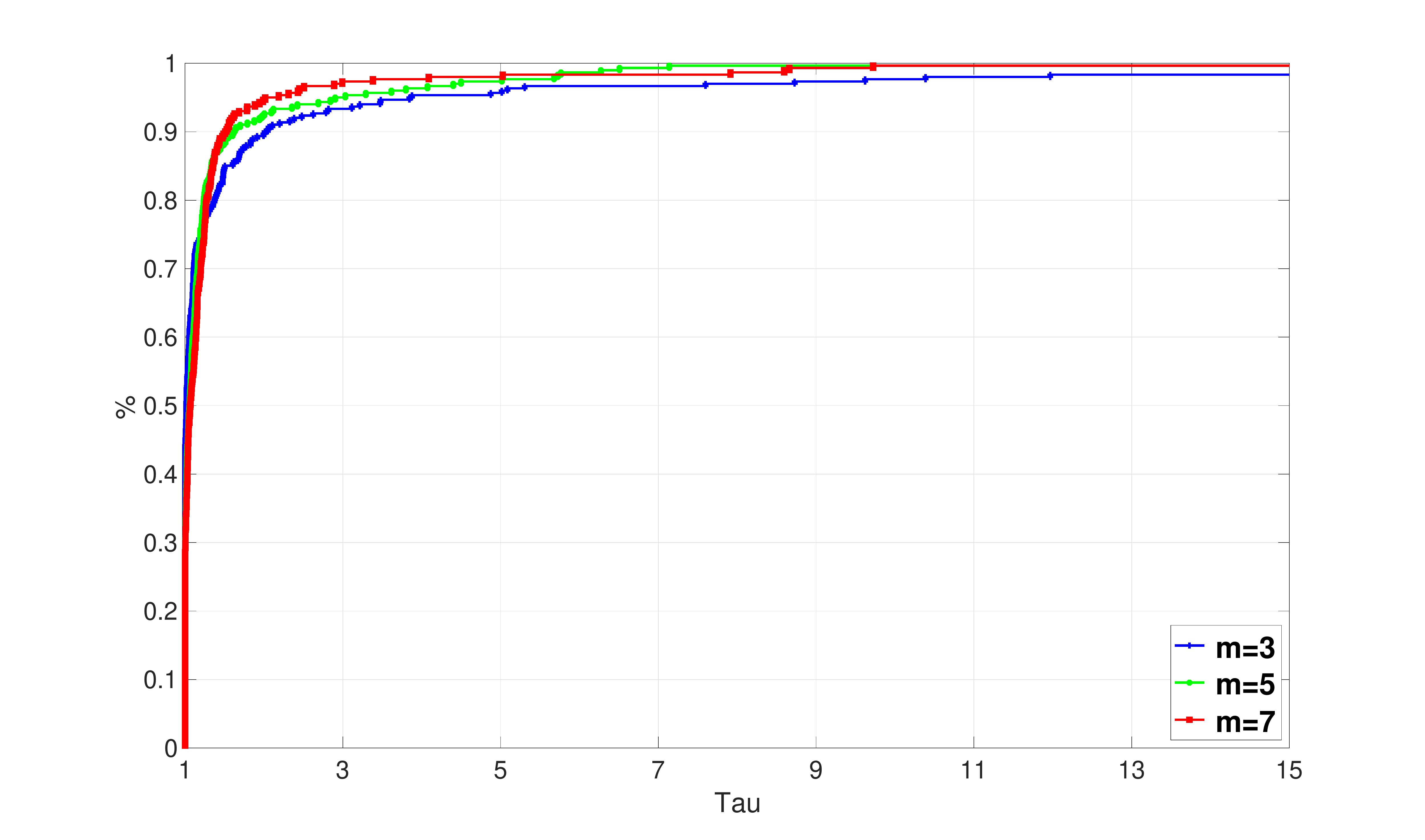}
\caption{Comparison of $m=3,5,7$.}
\label{figmemory}
\end{minipage}
\bigskip
\begin{minipage}{.5\textwidth}
\centering
\includegraphics[scale=0.12]{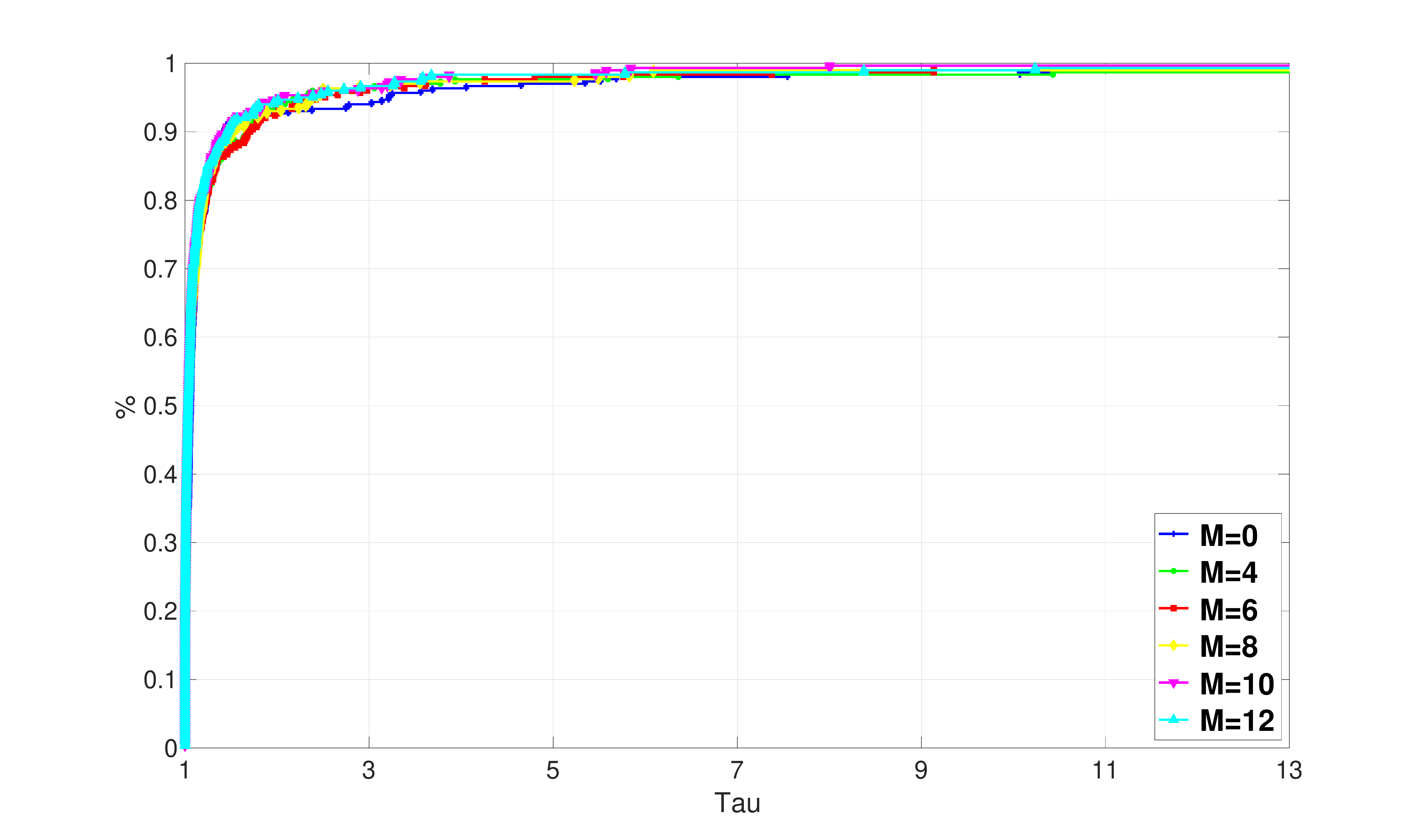}
\caption{Comparison of M.}\label{fignonmon}
\end{minipage}
\end{figure}

First we consider $\gamma_1$ and $\gamma_2$ that control regularized parameters.
We perform numerical experiments with $9$ different sets  of $(\gamma_1,\gamma_2)$ in Table \ref{tabgama}.
The remaining parameters are set to
$$ \eta_1=0.01, \eta_2=0.9, \mu_{min}=1.0\times 10^{-3},m=5,
M=10.$$
 Table \ref{tabgama} shows the 
number of success and rate of success for all $313$ problems. Figure \ref{figgama} shows the 
distribution function of these parameter sets in terms of the CPU time.

\begin{table}[H]
\caption{The number of success and rate of success at each $(\gamma_1,\gamma_2)$.}
\centering
\begin{tabular}{|c||c|c|c|c|}
\hline
P & $\gamma_1$ & $\gamma_2$ & Number of successes  & Success rate~(\%) \\
\hline
\hline
$P_1$ & $0.1$ & $2.0 $ & $261 $ & $ 83.4$ \\
\hline
$P_2$ & $0.1$ & $5.0 $ & $ 263$ & $ 84$ \\
\hline
$P_3$ & $0.1$ & $10.0 $ & $262 $ & $ 83.7$ \\
\hline
$P_4$ & $0.2$ & $2.0 $ & $261 $ & $ 83.4$ \\
\hline
$P_5$ & $0.2$ & $5.0 $ & $258$ & $ 82.4$ \\
\hline
$P_6$ & $0.2$ & $10.0 $ & $260 $ & $ 83$ \\
\hline
$P_7$ & $0.5$ & $2.0 $ & $256 $ & $ 81.8$ \\
\hline
$P_8$ & $0.5$ & $5.0 $ & $261 $ & $ 83.4$ \\
\hline
$P_9$ & $0.5$ & $10.0 $ & $262 $ & $ 83.7$ \\
\hline
\end{tabular}
\label{tabgama}
\end{table}  
\noindent
From Table \ref{tabgama} and Figure \ref{figgama}
it is clear that $(\gamma_1,\gamma_2)=(0.1,10.0) $ is the best. Therefore, we set 
$\gamma_1=0.1$ and $\gamma_2=10.0$ for all further experiments.


Next, we compare the number $m$ of vector pairs in the L-BFGS procedure. Note that  
the original \lb usually choose it  in $3\leq m \leq 7$ \cite{Nocedal}. Thus, we compare 
$m=3,5, 7$. The remaining parameters  are set to

$$\gamma_1=0.1,\gamma_2=10, \eta_1=0.01,\eta_2=0.9,\mu_{min}=1.0\times 
10^{-3},M=10. 
$$

\begin{table}[H]
\caption{The number of success and rate of success at each $m$.}
\centering
\begin{tabular}{|c||c|c|}
\hline
\textbf{Memory}  & Number of successes  & Success rate~(\%) \\
\hline
\hline
3   & 257   & 82.1\\
\hline
5 &   260  & 83\\
\hline
7 & 267 &85.3\\
\hline
\end{tabular}
\label{tabmemory}
\end{table}

From Table \ref{tabmemory} we see that $m=7$ is the best, while Figure \ref{figmemory} 
shows that $m=5$ is initially better in terms of CPU time. Therefore, we set
$m=5$ for further experiments.

Finally, we compare the behavior of nonmonotone parameters $M$. 
We compare  $M=0,4,6,8,10,12$. Note that $M=0$ implies the usual monotone decreasing
case. The remaining parameters are set to
$$\gamma_1=0.1,\gamma_2=10, \eta_1=0.01,\eta_2=0.9,\mu_{min}=1.0\times 
10^{-3},m=5.$$
Figure \ref{fignonmon} shows the distribution function of the nonmonotone parameter in terms of the CPU time.

From Table \ref{tabnonmon} and Figure \ref{fignonmon} it is clear that $M=10 $ is better.  Therefore, we use 
$M=10 $ in the next section.
\begin{table}[H]
\caption{The number of success and rate of success at each $M$.}
\centering
\begin{tabular}{|c||c|c|}
\hline
\textbf{Nonmonotone}  & Number of successes  & Success rate~(\%) \\
\hline
\hline
Monotone(M=0)   & 263    & 84\\
\hline
4 &   260  & 83\\
\hline
6 & 262 & 83.7\\
\hline
8 & 260 & 83\\
\hline
10 & 263 & 84\\
\hline
12 & 263 & 84\\
\hline
\end{tabular}
\label{tabnonmon}
\end{table}

\subsection[Comparisons of RL-BFGS-SW, RL-BFGS and L-BFGS method]{Comparisons of RL-BFGS-SW, RL-BFGS and L-BFGS \\method}

We compare the RL-BFGS-SW, the RL-BFGS and the L-BFGS methods in terms of function 
evaluations and CPU time. For all  numerical results, the parameters in RL-BFGS and
 RL-BFGS-SW are as follows:
$$\eta_1=0.01, \eta_2=0.9, 
\mu_{min} = 1.0 \times 10^{-3}, M=8,m=5, \gamma_1=0.1,  \gamma_2=10.$$

\begin{figure}[H]
\centering
\begin{minipage}[b]{.48\textwidth}
\centering
\includegraphics[scale=0.12]{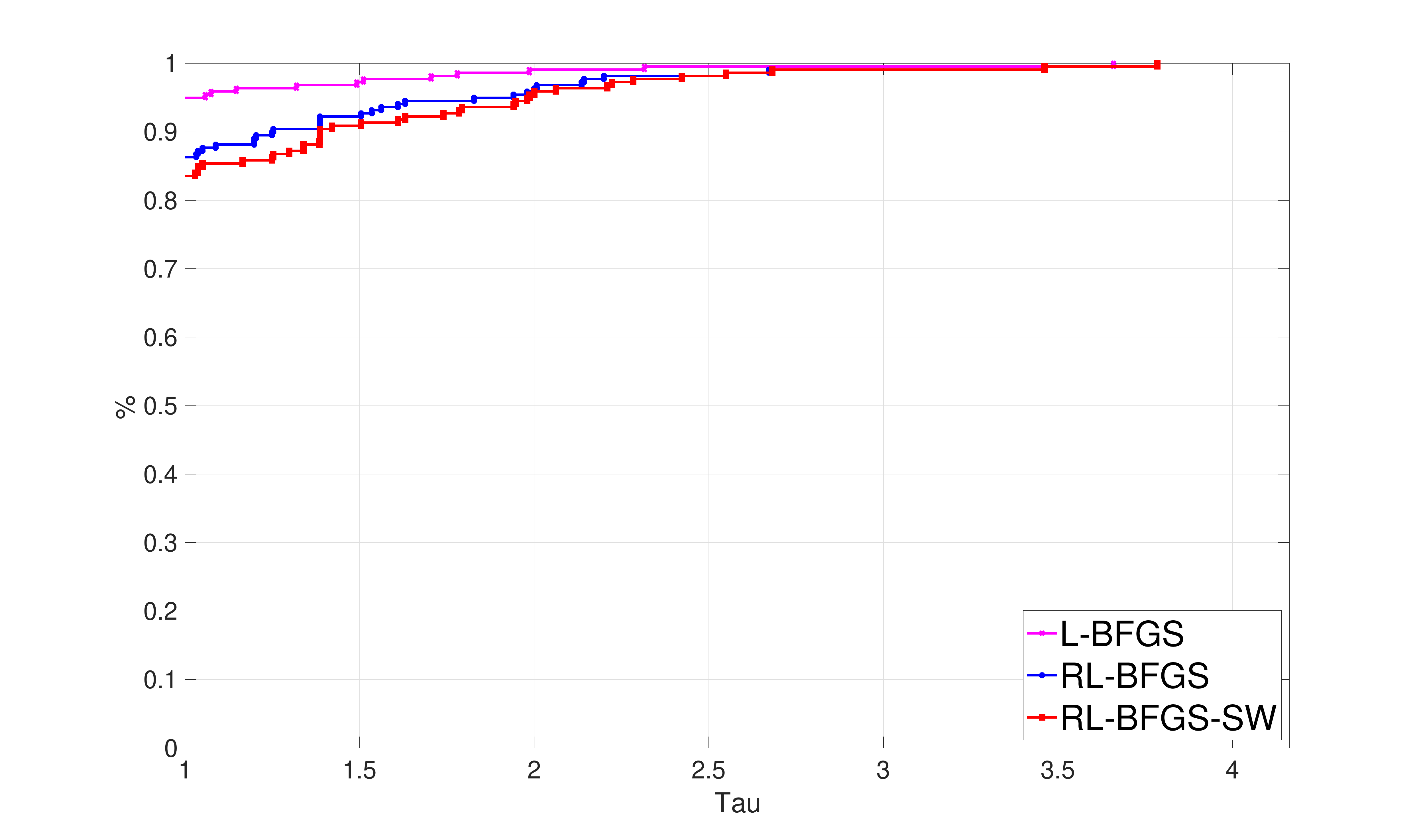}
\caption{Comparison of $n_f$.}
    \label{LBFGS1}
\end{minipage}%
\begin{minipage}[b]{.48\textwidth}
\centering
\includegraphics[scale=.12]{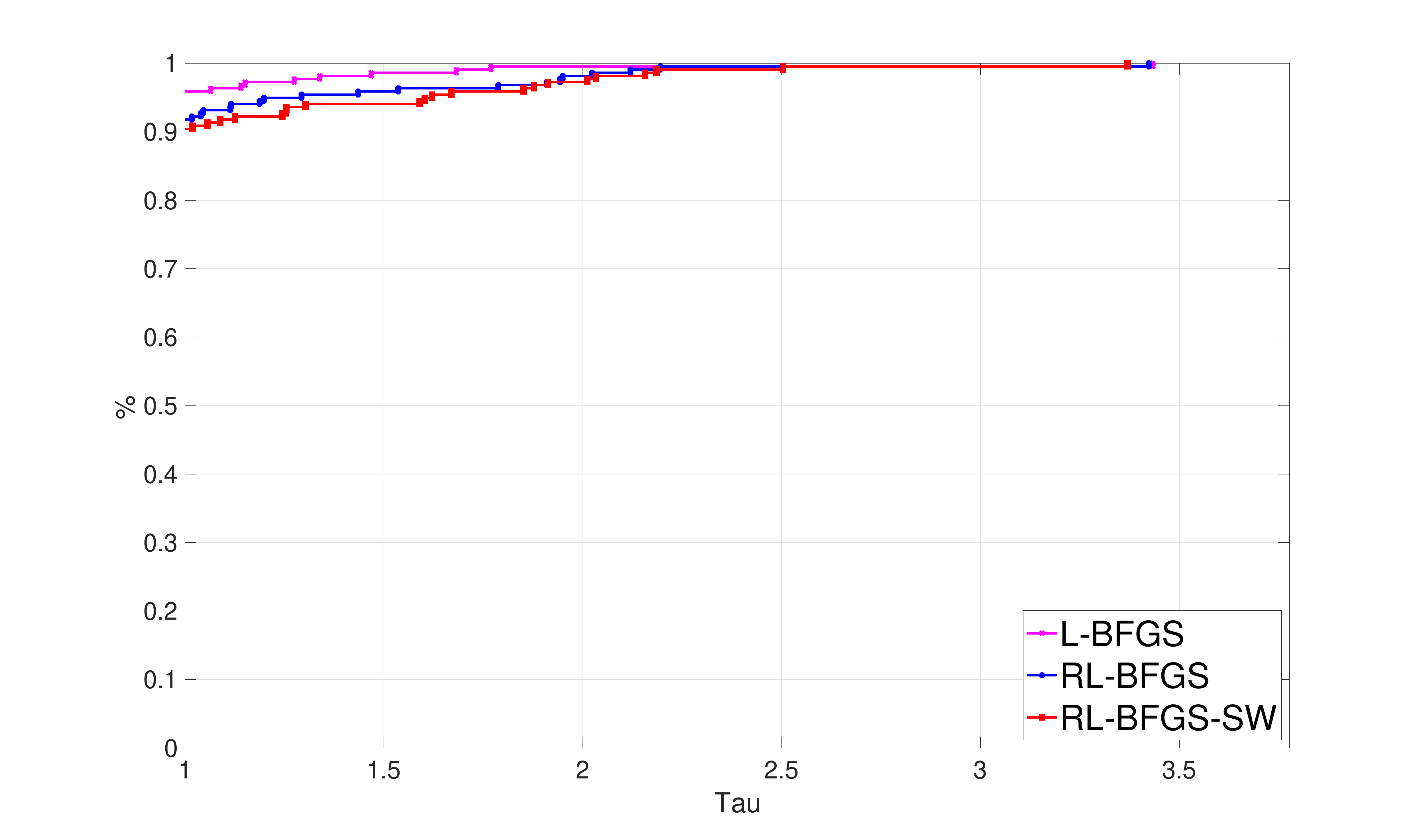}
\caption{Comparison  of CPU time.}
\label{LBFGS2}
\end{minipage}
\end{figure}

Table \ref{Tablen2n1} shows the results of the number of successes and rate of successes for 
all 313 test problems. Figures \ref{LBFGS1} and \ref{LBFGS2} show the results of
$\mathcal{P}_{\mathcal{S}}$ in terms of function evaluations and CPU time, respectively. Here 
$\mathcal{S}$ is the set of problems that are solved by all three algorithms. 

\begin{table}[H]
\caption{The number of success and rate of success for 313 problems.}
\centering
\begin{tabular}{|c||c|c|}
\hline
\textbf{Algorithm}  & Number of successes  & Success rate~(\%) \\
\hline
\hline
L-BFGS   & 225    & 71.9\\
\hline
RL-BFGS &   261  & 83.4\\
\hline
RL-BFGS-SW & 261 & 83.4\\
\hline
\end{tabular}
\label{Tablen2n1}
\end{table}

Table \ref{Tablen2n1} shows  that L-BFGS can solve 71.9\%  of test problems while 
both RL-BFGS and RL-BFGS-SW can solve 83.4\% of problems. On the other hand, Figures 
\ref{LBFGS1} and \ref{LBFGS2} show that L-BFGS is faster than the regularized ones for the solved problems.

\begin{figure}[H]
\begin{minipage}[b]{.48\textwidth}
\centering
\includegraphics[scale=0.12]{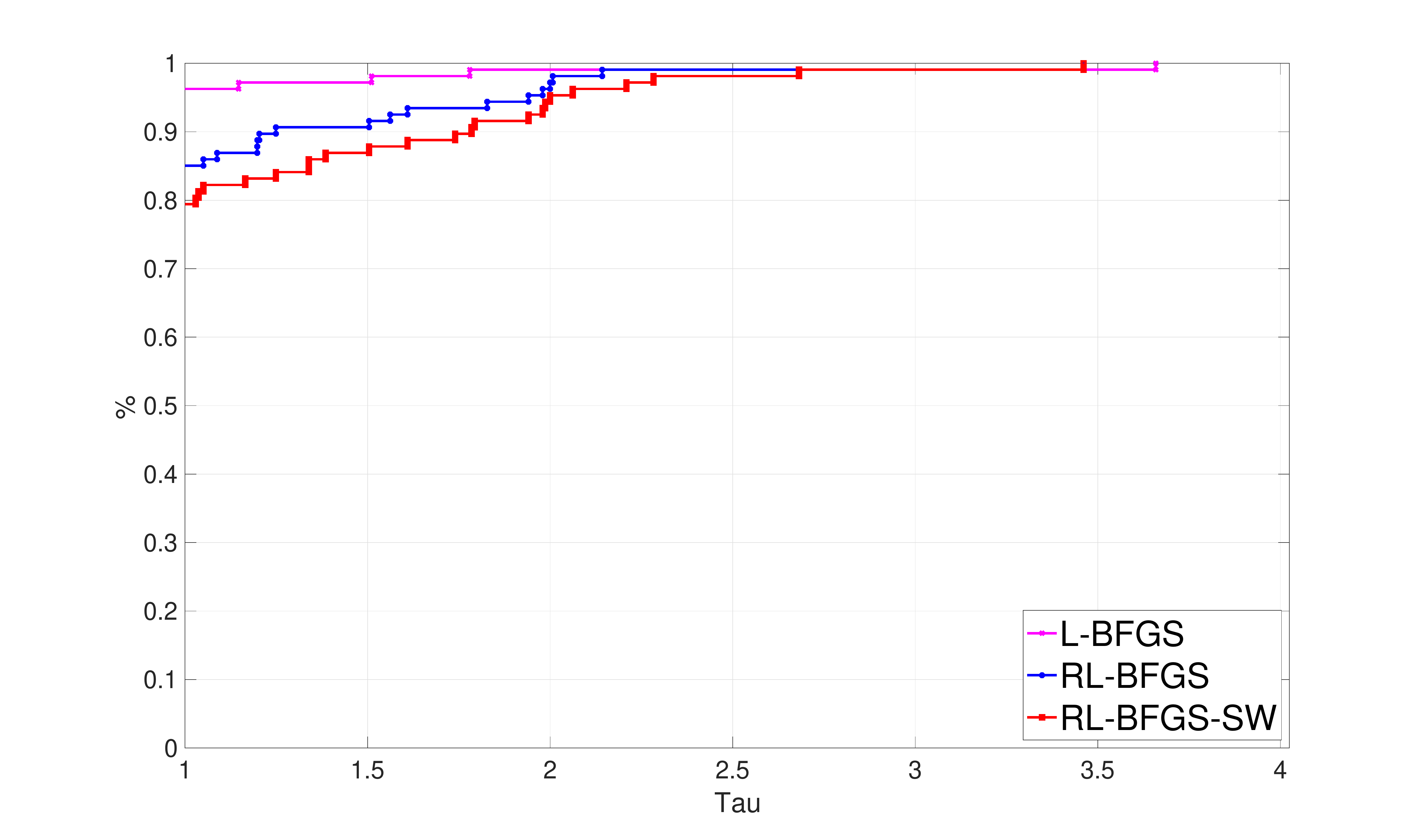}
\caption{Comparison  of $n_f$~(LS).}
    \label{LBFGS3}
\end{minipage}%
\begin{minipage}[b]{.48\textwidth}
\centering
\includegraphics[scale=.12]{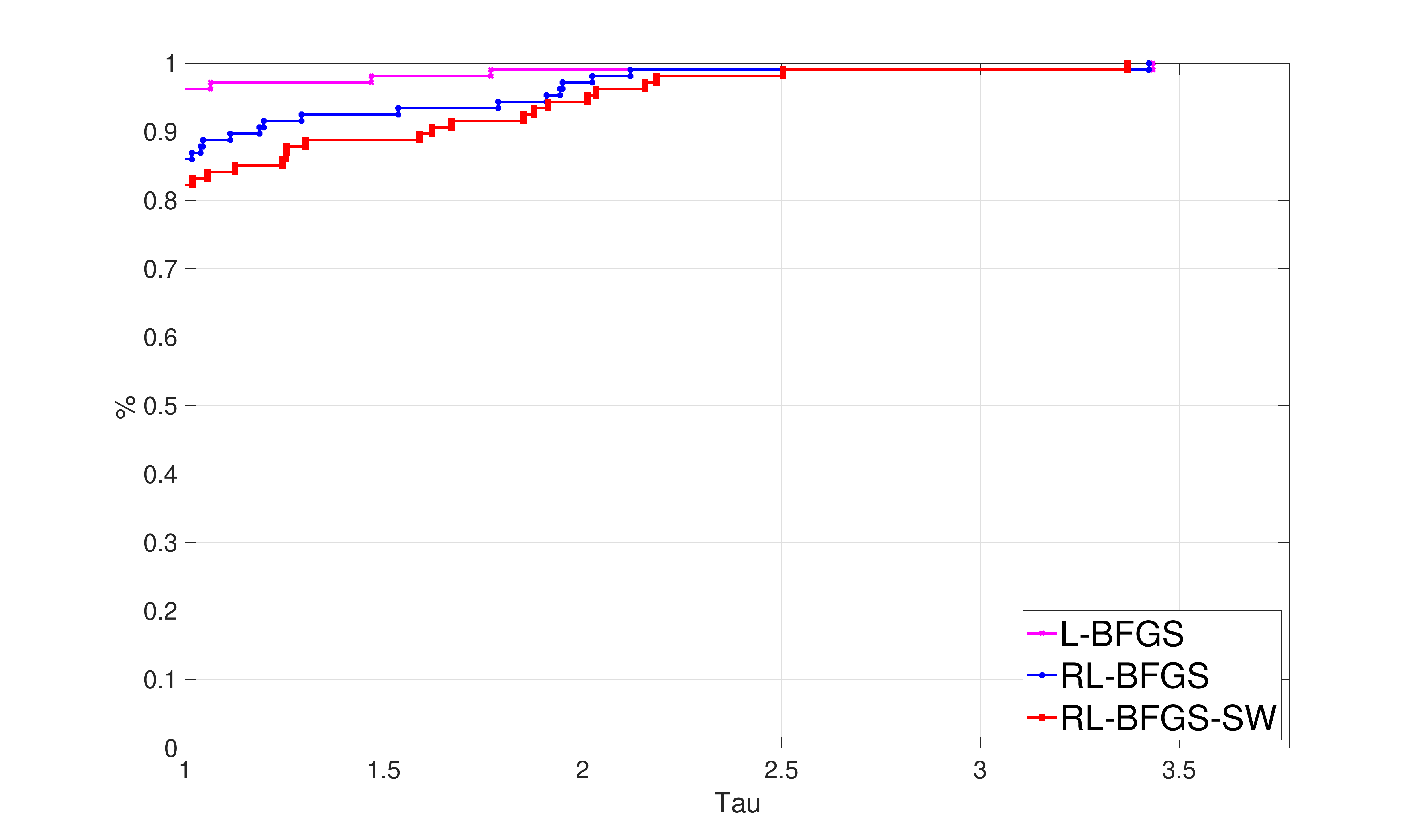}
   \caption{Comparison of CPU time~(LS).}
    \label{LBFGS4}
\end{minipage}
\end{figure}

We define the large-scale problem whose dimension is over or equal to $1000$. Table
\ref{Tablen23} shows the number of success and rate of success for the 151 large-scale 
problems. Furthermore, Figures \ref{LBFGS3} and \ref{LBFGS4} shows performances for 
$\mathcal{P}_{\mathcal{S}^{large}}$, where $\mathcal{P}_{\mathcal{S}^{large}}$ denotes all the
151 large-scale test problems from the $\mathcal{P}_{\mathcal{S}}$.

\begin{table}[H]
\caption{The number of success and rate of success for 151 large-scale problems.}
\centering
 \begin{tabular}{|c||c|c|}
\hline
\textbf{Algorithm}  & Number of successes  & Success rate~(\%) \\
\hline
L-BFGS   & 107  & 70.9\\
\hline
RL-BFGS &   125  & 82.8\\
\hline
RL-BFGS-SW & 125 & 82.8\\
\hline
\end{tabular}
\label{Tablen23}
\end{table}

Table \ref{Tablen23}  shows that the L-BFGS can solve 70.9\%  of test problems while both 
the RL-BFGS and RL-BFGS-SW can solve
82.8\%  of test problems. It concludes that both proposed 
methods can solve more number  of test problems as compare to the L-BFGS. On the other 
hand the above figures show 
that the L-BFGS requires fewer function evaluations than the proposed 
method. 

The above numerical results indicate that the numerical behaviors of the RL-BFGS and
the RL-BFGS-SW are almost same. 
To see the differences we present the numerical results which compare the performance of each test problem.

We observed that the RL-BFGS-SW performs the line search for 93 problems, and does not use it for the remaining
problems. Therefore, we compare the results for those 93 problems. Table \ref{ovstwlf12} shows the comparison in
terms of the number of function evaluations and Algorithm X~$<$~Algorithm Y means that the number of function
evaluations of the Algorithm X is fewer than that of the Algorithm Y. From the Table \ref{ovstwlf12}  we see 
that RL-BFGS-SW requires fewer number of function evaluations than that of RL-BFGS for 38 test problems while 
RL-BFGS requires fewer number of function evaluations than that of RL-BFGS-SW for 34 test problems among 93 test 
problems. Moreover, for the large-scale test problems, RL-BFGS-SW requires fewer number of function evaluations 
than that of RL-BFGS for 13 test problems while RL-BFGS requires fewer number of the function evaluations than 
that of RL-BFGS-SW for 20 test problems among 39 large-scale test problems. It concludes 
that the RL-BFGS with line search works well for some problems.



\begin{table}[H]
\caption{comparison for 93  problems solved by at least one algorithm
 (RL-BFGS or RL-BFGS-SW) in terms of $n_f$.}
\centering
 \begin{tabular}{|c||c|c|}
\hline
Number of problems  &  RL-BFGS-SW$<$RL-BFGS & RL-BFGS$<$RL-BFGS-SW \\
\hline
\hline
93 &   38  & 34\\
\hline
39(large-scale) & 13 & 20 \\
\hline
\end{tabular}
\label{ovstwlf12}
\end{table}
\section{Conclusion}
In this paper we have proposed a combination of the \lb and the regularization technique. We showed the global
convergence under appropriate assumptions. We have also presented some efficient implementations. In numerical
results, the overall comparison shows that the proposed method can solve more problems than the original L-BFGS. 
This result indicates that the proposed method is robust in terms of solving number of 
problems. 

For future work, we may consider proposing the stochastic version of the proposed method to solve empirical risk
minimization problems.


\end{document}